\documentclass[final,5p,times,twocolumn,authoryear]{elsarticle}

\usepackage{amssymb}
\usepackage{array}
\usepackage{amsmath}
\usepackage{amsfonts}
\usepackage{adjustbox}
\usepackage{float}
\usepackage{caption}
\usepackage{graphicx}
\usepackage{subcaption}
\usepackage{booktabs}
\usepackage{makecell}
\usepackage{bm}
\usepackage{multirow}
\usepackage{enumitem}
\usepackage[colorlinks,
            linkcolor=red,
            anchorcolor=blue,
            citecolor=green
            ]{hyperref}
\allowdisplaybreaks

\newcommand{\picref}[1]{Fig.\ \ref{#1}}
\newcommand{\secref}[1]{Section\ \ref{#1}}
\newcommand{\tabref}[1]{Table\ \ref{#1}}

\journal{Computers \& Operations Research}

\begin{document} 

\begin{frontmatter}


\title{Towards An Unsupervised Learning Scheme for Efficiently Solving Parameterized Mixed-Integer Programs}

\author[THU]{Shiyuan Qu}
\author[PPEI,CNPC]{Fenglian Dong}
\author[PPEI,CNPC]{Zhiwei Wei}
\author[THU]{Chao Shang}

\affiliation[THU]{organization={Department of Automation, Beijing National Research Center for Information Science and Technology, Tsinghua University},
            city={Beijing},
            postcode={100084}, 
            country={China}}

\affiliation[PPEI]{organization={PetroChina Planning and Engineering Institute},
            city={Beijing},
            postcode={100083}, 
            country={China}}

\affiliation[CNPC]{organization={CNPC Laboratory of Oil \& Gas Business Chain Optimization},
            city={Beijing},
            postcode={100086}, 
            country={China}}

\begin{abstract}
In this paper, we describe a novel unsupervised learning scheme for accelerating the solution of a family of mixed integer programming (MIP) problems. Distinct substantially from existing learning-to-optimize methods, our proposal seeks to train an autoencoder (AE) for binary variables in an unsupervised learning fashion, using data of optimal solutions to historical instances for a parametric family of MIPs.
By a deliberate design of AE architecture and exploitation of its statistical implication, we present a simple and straightforward strategy to construct a class of cutting plane constraints from the decoder parameters of an offline-trained AE. These constraints reliably enclose the optimal binary solutions of new problem instances thanks to the representation strength of the AE. More importantly, their integration into the primal MIP problem leads to a tightened MIP with the reduced feasible region, which can be resolved at decision time using off-the-shelf solvers with much higher efficiency. 
Our method is applied to a benchmark batch process scheduling problem formulated as a mixed integer linear programming (MILP) problem.
Comprehensive results demonstrate that our approach significantly reduces the computational cost of off-the-shelf MILP solvers while retaining a high solution quality. 
The codes of this work are open-sourced at \href{https://github.com/qushiyuan/AE4BV}{https://github.com/qushiyuan/AE4BV}.
 
\end{abstract}

\begin{keyword}
Mixed integer programming \sep Machine learning \sep Autoencoder \sep Scheduling 
\end{keyword}

\end{frontmatter}

\section{Introduction} 

Mixed integer programming (MIP) problems comprise a significant part of mathematical programming, in which a portion of decision variables are constrained to take integer values. 
Due to their remarkable modeling power, MIPs have served as useful tools in diverse fields such as 
scheduling \citep{sawik2011scheduling, ku2016mixed}, planning \citep{pochet2006production, HASSAN2021105174}, routing \citep{ vidovic2014mixed, WANG2025106883}, etc. 
Although binary variables enhance the descriptive and representational power, they also substantially complicate the solution process of MIP models \citep{wolsey2007mixed}. 
When both the objective and constraint functions of MIP are linear, the MIP problem becomes a Mixed-Integer Linear Programming (MILP) problem.  
Currently, off-the-shelf commercial and open-source solvers are available for MILPs, such as Gurobi \citep{gurobi}, SCIP \citep{BolusaniEtal2024OO} , and COPT \citep{copt}, which primarily build upon the branch-and-cut (B\&C) algorithm and offer certificates of optimality. Despite these successes, MIPs remain a class of optimization problems that are challenging to solve, especially for large-scale instances.

Recent years have been witnessing a surging interest in \textit{learning to optimize} (L2O), which leverages machine learning tools to speed up the solution of complicated optimization problems \citep{bengio2021machine, zhang2023survey}. 
An important line of work has focused on accelerating the B\&C process of MILPs with machine learning-based branching rules that improve upon traditional heuristics.
The earliest works in this vein focused on using machine learning to imitate established branching heuristics, with the goal of reducing computational costs and enhancing performance \citep{marcos2014supervised, gasse2019exact}. 
\cite{zarpellon2021parameterizing} proposed to utilize structural knowledge of the branching tree to improve the generalization capability of machine learning models. 
For efficient adaptation to out-of-sample problem instances, online learning of branching heuristics has also been explored by \cite{khalil2016learning}. 
Additionally, similar ideas have also been adopted in learning data-based heuristics for cut selection within the B\&C process \citep{cornuejols2008valid, deza2023machine}. 
For instance, \cite{tang2020reinforcement, paulus2022learning} investigated methods for scoring individual cuts to enhance cut selection efficiency, and this idea was further extended to scenarios involving multi-type cuts \citep{huang2022learning}. 
For a better trade-off between the quality and the cost of cut selection, a regression forest approach was suggested in \cite{berthold2022learning}. 

The strength of machine learning can also be leveraged by predicting optimal values of binary variables in MIPs in an end-to-end fashion. 
The key idea is to build a mapping from problem parameters to optimal solutions by performing supervised learning on a labeled dataset. 
In this way, by fixing binary variables to their predictions, the remaining problem becomes a much simpler (though possibly nonconvex) program involving continuous variables only, and thus it is no longer necessary to solve the entire MIP from scratch.
In \cite{li2018combinatorial}, 
they proposed to reformulate MILPs into graphs and train graph convolutional networks to leverage structural knowledge for predicting and refining the solutions. 
For improved interpretability, \cite{bertsimas2021voice, bertsimas2023global} presented a novel end-to-end approach based on the classification tree and decision tree. 
However, these prediction methods may be prone to loss of optimality or even result in infeasible solutions. 
To mitigate this deficiency, a remedy is to use the predicted solution as a warm start for off-the-shelf MILP solvers; see e.g. \cite{xavier2021learning, sugishita2024use}. 

In many practical scenarios, one may face a multitude of optimization problems that involve binary variables and need to be solved repeatedly. 
Typical examples include security-constrained unit commitment (SCUC) problems \citep{xavier2021learning}, hybrid control problems \citep{dragotto2023differentiable}, etc.  
These instances typically have a common similar structure and thus can be viewed as a \textit{parametric family} of MIPs. 
As a consequence, their optimum binary variables can exhibit highly complex and implicit correlations.   
For instance, in scheduling problems cast as MIPs, binary variables have real-world physical implications such as task sequencing or resource allocation, and are optimized to reach the best possible solution. 
Driven by a common objective of cost minimization or efficiency maximization, the optimum for binary variables in a family of MIPs is unlikely to be evenly distributed across the entire feasible region but embodies a particular statistical pattern, exhibiting a ``data-rich, information-poor" characteristic. 
In other words, the optimal binary solutions tend to reside in a low-dimensional manifold rather than spreading the entire high-dimensional space of decision variables. 
To the best of our knowledge, such a latent structure has not been exploited yet  within the current L2O paradigm. 

To fill this gap, we seek to disentangle the latent structure of optimal binary solutions of a parametric family of MIPs in order to effectively expedite the solution of new instances using MIP solvers. Following the spirit of unsupervised learning, we propose to train an auto-encoder (AE) offline for efficient dimension reduction of optimal binary variables for a family of MIPs. 
By a deliberate design of AE architecture and exploitation of its generative nature, a simple and straightforward approach is proposed to construct a class of cutting planes from network parameters of AE, in the form of a polytopic set that reliably encloses the optimal binary solutions of new problem instances due to the strength of autoencoder. 
By adding these cutting planes to the primal problem, tightened MIPs can be formulated for new target instances and tackled using off-the-shelf MIP solvers at a much higher speed. 
Specifically, thanks to the reduced feasible region of binary variables, the convergence of B\&C process of MILP can be considerably accelerated and thus the expensive computational cost gets effectively reduced. 

In a nutshell, our method offers a new and effective route towards solving complicated MIPs when existing end-to-end learning methods struggle with ensuring feasibility and off-the-shelf solvers are still necessary. Different from current learning-based heuristics for branching and warm starting, our method lies at the unexplored intersection between unsupervised learning and global optimization, thereby broadening the scope of L2O. Its usefulness is comprehensively evaluated on a batch process scheduling problem, which can be recast as a complicated MILP and handled using off-the-shelf solvers. In this application, existing end-to-end L2O methods are shown to fall short of producing feasible solutions. By contrast, our unsupervised learning approach not only ensures feasibility but also achieves an average improvement of 90\% in the solution speed, while maintaining high solution quality.

The remainder of this paper is structured as follows.
In Section 2, we offer the problem formulation, introducing the basic idea and the motivation of our approach. 
In Section 3, we detail our method, explaining how to learn the latent structure of optimal solutions from historical instances using an AE and how to construct a class of cutting plane constraints from the decoder parameters in a well-trained AE. 
Section 4 evaluates our approach through applications to a batch processing problem.
We conclude in Section 5 by summarizing our contributions and discussing potential directions for future work. 

\section{Problem Formulation and Motivation}

The standard formulation of an MIP comprises a set of decision variables and an objective function while being subject to certain constraints. 
Formally, a family of MIPs \textit{parameterized by} $\boldsymbol{\theta}$ can be expressed as:
\begin{equation}
    \begin{aligned}
          {\rm MIP}[\boldsymbol{\theta}] : & &\min_{\boldsymbol{x}, \boldsymbol{u}}~~ & g(\boldsymbol{x}, \boldsymbol{u};\boldsymbol{\theta})\\
        && \mbox{s.t.}~~~ &  f(\boldsymbol{x}, \boldsymbol{u}; \boldsymbol{\theta}) \geq \boldsymbol{0}\\
        & &&  \boldsymbol{x} \in \mathbb{R} ^ {m}, \boldsymbol{u}\in \mathbb{B} ^{p}\\
    \end{aligned} 
    \label{standard mip}
\end{equation}
where $g(\cdot)$ denotes the objective function and $f(\cdot)$ denotes the constraint functions, $\boldsymbol{x}$ is the vector of continuous decision variables, and $\boldsymbol{u}$ is the vector of binary decision variables. $\boldsymbol{\theta}$ denotes problem parameters that may vary across different instances, giving rise to a family of similar MIPs. 
When $g(\cdot)$ and $f(\cdot)$ are linear functions, the MIP[$\boldsymbol{\theta}$] becomes an MILP and can represent a wide spectrum of decision-making problems, including SCUC problems \citep{xavier2021learning}, supply chain planning problems \citep{misra2022learning}, and voltage regulation problems \citep{sang2023online}. 

Data collected from solving multiple instances of MIPs contains information that is useful for solving new problem instances. Assume that there are $N$ MIPs, i.e. MIP[$\boldsymbol{\theta}_{n}$], $n=1,\cdots,N$, which are parameterized by $\boldsymbol{\theta}_n,~n=1,\cdots,N$ and have already been solved to optimality, yielding a pre-collected dataset of $N$ optimal solutions $(\boldsymbol{x}^*_n, \boldsymbol{u}^*_n),~n=1,\cdots,N$. 
The binary variables $\boldsymbol{u}$ are the main source of computational complexity in solving MIPs \citep{conejo2006decomposition}. 
Once these binary variables are fixed, the resulting problem only involves continuous variables and thus becomes much easier to solve. Consequently, a natural idea is to build an end-to-end mapping $\mathcal{F}(\cdot)$ from $\boldsymbol{\theta}$ to $\boldsymbol{u}^*$ by means of supervised learning on a labeled dataset  $\mathcal{D} = \{(\boldsymbol{\theta}_n, \boldsymbol{u}^*_n),~n=1,\cdots,N\}$. This mapping can be used to infer optimal binary variables $\boldsymbol{u}^*_{\rm new}$ for a new problem instance MIP[$\boldsymbol{\theta}_{\rm new}$], without the need for off-the-shelf MIP solvers at decision time. Despite the significant reduction in solution time, these end-to-end L2O methods require a huge volume of labeled data. Besides, it is not trivial to ensure the feasibility of the predicted binary variables, especially when tackling large-scale MIPs with complicated constraints. 

In this article, we present a novel L2O approach harnessing the strength of unsupervised learning to disentangle latent structures in the optimum of binary variables. 
This is formalized as training an AE from an unlabeled dataset $\mathcal{D}_{\boldsymbol{u}} = \{\boldsymbol{u}^*_n,~n=1,\cdots,N\}$ for dimension reduction and feature extraction of binary variables $\boldsymbol{u}^*$. This can be effectively done in an offline manner, due to the rapid advancement of training deep neural networks. 
Interestingly, we show that one can straightforwardly construct \textit{cutting planes}, a fundamental concept for solving MIPs since the 1990s \citep{contardo2023cutting}, from the decoder parameters of a well-trained AE. 
These data-based cutting planes serve as polytopic constraints that enclose the optimal solutions of a parametric family of MIPs with high probability. 
By embedding these polytopic constraints into the primal problem, we arrive at a \textit{tightened MIP}: 
\begin{subequations}
    \begin{align}
        {\rm T\text{-}MIP}[\boldsymbol{\theta}] : &&\min_{\boldsymbol{x}, \boldsymbol{u}}~~ & g(\boldsymbol{x}, \boldsymbol{u};\boldsymbol{\theta})& \\
        && \mbox{s.t.}~~~ &f(\boldsymbol{x}, \boldsymbol{u}; \boldsymbol{\theta}) \geq \boldsymbol{0}& \\
        && & \boldsymbol{A}\boldsymbol{u} + \boldsymbol{b} \geq \boldsymbol{0} \label{cutting planes in standard MIP} & \\
        && & \boldsymbol{x} \in \mathbb{R} ^ {m}, \boldsymbol{u}\in \mathbb{B} ^{p} & 
    \end{align}
    \label{standard MIP with cutting planes}
\end{subequations}
When dealing with a new target problem with parameter $\boldsymbol{\theta}_\text{new}$, the tightened problem T-MIP[$\boldsymbol{\theta}_\text{new}$] can be solved in the same way but with higher efficiency than its primal version MIP[$\boldsymbol{\theta}_\text{new}$]. 
Importantly, the feasible region of binary variables in T-MIP[$\boldsymbol{\theta}_\text{new}$] can be effectively reduced thanks to the added cuts \eqref{cutting planes in standard MIP} that encode useful information of the distribution of $\boldsymbol{u}^*$.
In particular, when the primal problem is an MILP, the integration of these cutting planes enables faster convergence of B\&C process of off-the-shelf MILP solvers. 

\section{Learning cutting planes using AE}\label{sec: method}

This section describes our unsupervised learning scheme to derive cutting planes by training AE from pre-collected data. 

\subsection{Autoencoder for binary variables (AE4BV)} \label{sec: AE4BV}

For a parametric family of MIPs, the optimal values of binary variables are not evenly distributed across the entire feasible region. 
To address this, we propose to utilize an AE with binary inputs and outputs to capture the latent structure in the optimal solutions. 
As a powerful machine learning model, AE as well as its variants has been extensively utilized for dimension reduction and feature extraction \citep{hinton2006reducing, kingma2013auto}.
It has a well-known encoder-decoder architecture, in which the encoder reduces the dimensionality of the input vector $\boldsymbol{u} \in \mathbb{B} ^{p}$, effectively compressing it into a lower-dimensional feature representation $\boldsymbol{h} = [ h_1 ~ h_2 ~ \cdots ~ h_d ]^\top \in \mathbb{R} ^{d}$ with $d < p$. 
The decoder reconstructs the output vector $\boldsymbol{v}$ from the feature $\boldsymbol{h}$, minimizing the reconstruction error between $\boldsymbol{u}$ and $\boldsymbol{v}$. 
The intermediate hidden layer produces the latent feature $\boldsymbol{h}$ as an \textit{information bottleneck}, forcing the network to learn representations that are compact but necessary for reconstructing the inputs. 

\begin{figure}[!htbp]
    \centering
    \includegraphics[width=\linewidth]{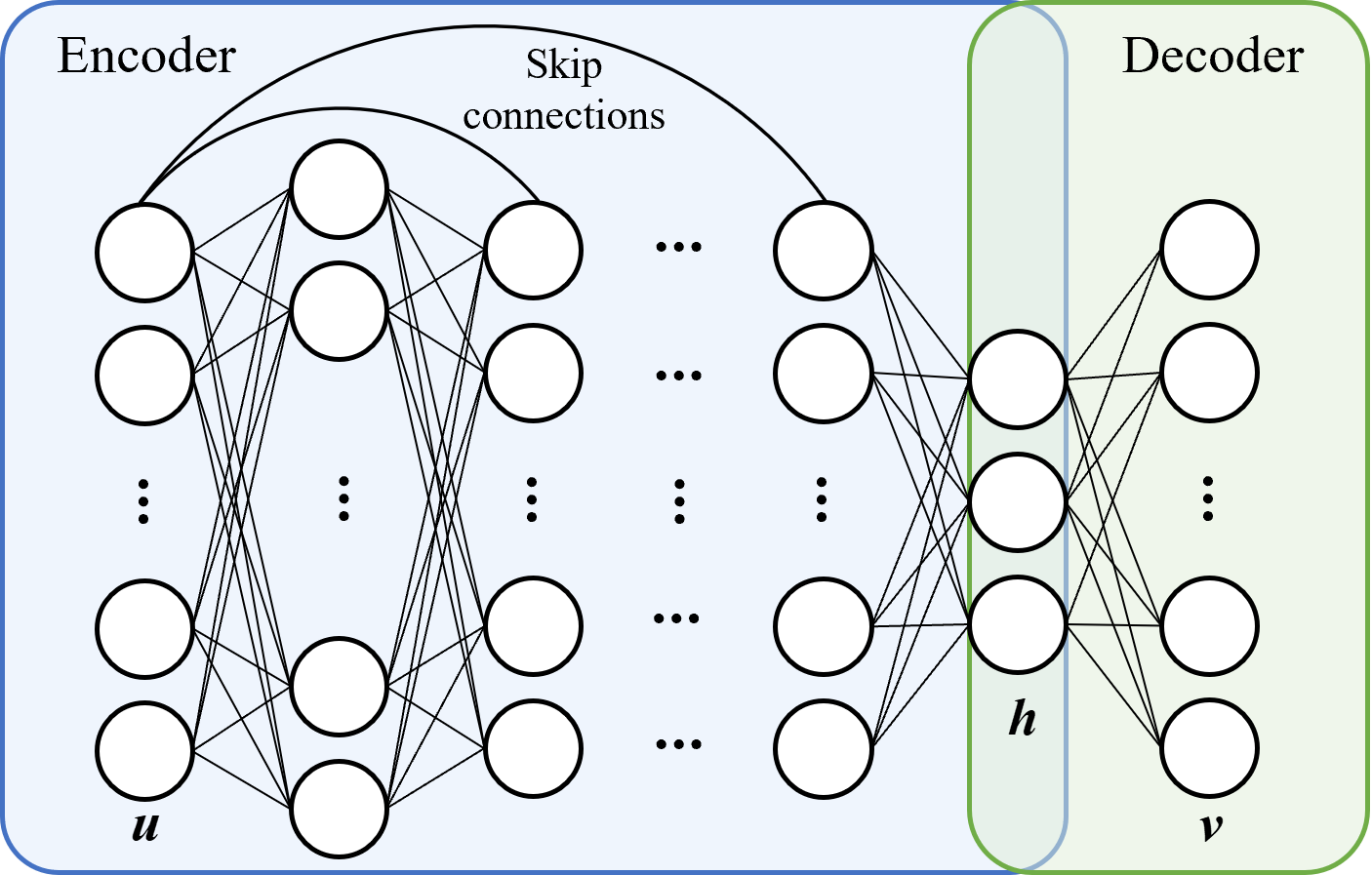}
    \caption{The network structure of AE4BV. }\label{fig:AE4BV}
\end{figure}

In this work, we employ an autoencoder tailored to binary variables, referred to as AE4BV, to achieve efficient dimension reduction and feature extraction of binary data. The network architecture is visualized in 
\picref{fig:AE4BV}, where the encoder consists of multiple layers. 
To regulate both training and inference complexity, all layers within AE4BV are fully connected. 
Moreover, to enhance the network efficiency while maintaining its representation capability, the fully connected layers in the encoder are augmented with skip connections  \citep{he2016deep}.

In stark contrast to most AE architectures, the decoder of AE4BV has only a single fully connected layer. As will become clear in \secref{sec: cutting plane constraints}, this simple single-layer structure is key to the successful construction of cutting planes from network parameters. Mathematically, the decoder can be expressed as:
\begin{equation}
    \boldsymbol{v} = \sigma(\boldsymbol{W}^\top\boldsymbol{h} + \boldsymbol{a})
    \label{decoder equation}
\end{equation}
where $\sigma(x) = 1/(1 + e^{-x}) \in (0,1)$ is the sigmoid function and is applied element-wise in \eqref{decoder equation}.
$\boldsymbol{W}$ and $\boldsymbol{a}$ denote the weight matrix and bias vector, respectively, with appropriate dimensions. 
It is noteworthy that $d$-dimensional feature $\boldsymbol{h}$ consists of continuous variables rather than binary variables. This is useful for better ability of dimension reduction and feature representation. 

The sigmoid function provides a probabilistic interpretation of outputs, which have been used in various models, such as logistic regression and restricted Boltzmann machines; see e.g. \citep{niculescu2005predicting, zou2019logistic, fischer2012introduction}. 
Given binary inputs $\boldsymbol{u} \in \mathbb{B} ^{p}$, we aim to reconstruct the inputs $\hat{\boldsymbol{u}} = [ \hat{u}_1 ~ \hat{u}_2 ~ \cdots ~ \hat{u}_p ]^\top$ that are binary as well. 
From a probabilistic perspective, the output \( v_i \) of the sigmoid function lies within the interval $(0,1)$, which can be interpreted as the probability of seeing $\hat{u}_i = 1,~i = 1,\dots,p$. 
Consequently, the binary inputs $\boldsymbol{u}$ can be reconstructed from $\boldsymbol{v}$ according to the following principle:
\begin{equation}
    \hat{u_i} = \left\{
        \begin{aligned}
            &0, \ 0 < v_i \le 1/2\\
            &1, \ 1/2 < v_i < 1, ~i=1, \dots, p
        \end{aligned} \label{binarization equation}
    \right.
\end{equation}
This implies that the outputs of AE4BV can be interpreted as a \textit{multi-label classifier}. Thus, we adopt the cross-entropy (CE), a prevalent cost function in multi-label classification (MLC) \citep{mannor2005cross}, to evaluate the reconstruction performance \citep{liu2017learning}:  
\begin{equation}
    \mathcal{L}(\boldsymbol{u}, \boldsymbol{v}|\Theta) = -\sum_{i = 1}^P[v_i\cdot \log u_i + (1 - v_i) \cdot \log (1 - u_i)],
    \label{eq: MLC}
\end{equation} 
where $\Theta$ denotes network parameters to be optimized over, including all encoder parameters and decoder parameters $\{\boldsymbol{W}, \boldsymbol{a}\}$. Then AE4BV can be trained by minimizing the CE loss \eqref{eq: MLC} on a pre-collected dataset using modern deep learning toolboxes such as PyTorch \citep{Ansel_PyTorch_2_Faster_2024}. 

\subsection{Cutting plane generation via AE4BV}\label{sec: cutting plane constraints}
Upon training AE4BV from data, our attention is then placed on the decoder, which itself effectively establishes a compact mapping from the feature representation $\boldsymbol{h}$ to the output $\hat{\boldsymbol{u}}$. We seek to explicitly express the latent structure of binary inputs $\boldsymbol{u}$ in terms of a convex polytope. 
By combining \eqref{decoder equation} and \eqref{binarization equation}, 
the reconstruction principle underlying $\hat{\boldsymbol{u}}$ can be rewritten as: 
\begin{equation}
    \left\{
        \begin{aligned}
            &\boldsymbol{W}_i^\top\boldsymbol{h} + a_i \geq 0 \Rightarrow \hat{u}_i = 1\\
            &\boldsymbol{W}_i^\top\boldsymbol{h} + a_i \leq 0 \Rightarrow \hat{u}_i = 0, ~i=1, \dots, p\\
        \end{aligned}
    \right.
    \label{decoder piecewise inequality}
\end{equation}
where $\boldsymbol{W}_i^\top$ and $a_i$ represent the $i^{\rm th}$ row of $\boldsymbol{W}^\top$ and the $i^{\rm th}$ entry of $\boldsymbol{a}$, respectively.
For a well-trained AE4BV, the interpretations of the reconstruction principle \eqref{decoder piecewise inequality} are two-fold. 
On the one hand, given a binary input $\boldsymbol{u}$, \eqref{decoder piecewise inequality} describes the decoding principle in an encoding-decoding regime, where $\boldsymbol{h}$ is produced from the encoder. On the other hand, given an arbitrary feature $\boldsymbol{h}$, one can generate from \eqref{decoder piecewise inequality} a new sample of $\boldsymbol{u}$, which is expected to be governed by the same distribution as the training dataset $\mathcal{D}_{\boldsymbol{u}}$. This essentially follows a \textit{generative} perspective of employing auto-encoders, which has played an indispensable role in the well-known generative adversarial networks \citep{goodfellow2020generative}. Thus, \eqref{decoder piecewise inequality} can be interpreted as a set-based description of $\boldsymbol{u}$ by introducing $\boldsymbol{h}$ as auxiliary variables: 
\begin{equation}
    \mathcal{P}_{\boldsymbol{u}} = \left\{
    \boldsymbol{u}~\middle|~
    \begin{aligned}
    \exists \boldsymbol{h}~\text{s.t.}~ &\boldsymbol{W}_i^\top \boldsymbol{h} + a_i \geq 0 \Rightarrow u_i = 1, \\
    &\boldsymbol{W}_i^\top \boldsymbol{h} + a_i \leq 0 \Rightarrow u_i = 0,~i=1, \dots, p
    \end{aligned}
    \right\}.
    \label{decoder piecewise inequality in set}
\end{equation}

In order to recast \eqref{decoder piecewise inequality in set} into linear constraints, we apply the big-M strategy: 
\begin{equation}
    \begin{split}
    \mathcal{P}_{\boldsymbol{u}} &= \left\{
    \boldsymbol{u}~\middle|~
    \begin{aligned}
    \exists \boldsymbol{h}~\text{s.t.}~
    &\boldsymbol{W}_i^\top\boldsymbol{h} + a_i \geq M(u_i -1), \\
    &\boldsymbol{W}_i^\top\boldsymbol{h} + a_i \leq Mu_i,~i=1, \dots, p
    \end{aligned}
    \right\}\\
    &= 
    \left\{
    \boldsymbol{u}~\middle|~
    \begin{aligned}
    \exists \boldsymbol{h}~\text{s.t.}~
    & \boldsymbol{W}^\top\boldsymbol{h} + \boldsymbol{a} \geq M\cdot(\boldsymbol{u} - \mathbf{1}), \\
    &\boldsymbol{W}^\top\boldsymbol{h} + \boldsymbol{a} \leq M\cdot\boldsymbol{u}
    \end{aligned}
    \right\}\\
    \end{split}
    \label{decoder piecewise inequality with Big-M}
\end{equation}
where $M>0$ is a sufficiently large constant and $\mathbf{1}$ is an all-ones vector. 
The inequalities in \eqref{decoder piecewise inequality with Big-M} define a convex polytope $\mathcal{P}_{\boldsymbol{u}}$, which describes the distributional geometry of the optimal binary variables $\boldsymbol{u}^*$. 
Importantly, the left-hand coefficient in \eqref{decoder piecewise inequality with Big-M} is the decoder parameters $\{\boldsymbol{W}, \boldsymbol{a}\}$ which can be directly derived from a well-trained AE4BV. 
As for the selection of $M$, it can be guided by the empirical distribution $\mathcal{D}_{\boldsymbol{h}} = \{\boldsymbol{h}_n,~n=1,\cdots,N\}$, where $\boldsymbol{h}_n$ is the feature vector inferred from $\boldsymbol{u}_n^*$ based on the trained AE4BV. 
Then a simple heuristic is given by: 
\begin{equation}
\begin{aligned}
        M &= \max\left\{|\boldsymbol{W}_i^\top\boldsymbol{h}_n + a_i|,~i=1,\cdots,p, \boldsymbol{h}_n \in \mathcal{D}_h\right\}.
    \end{aligned}
    \label{eq: estimation M}
\end{equation}

\begin{figure}[htbp]
	\centering
	\begin{subfigure}{0.35\linewidth}
		\centering
		\includegraphics[width=0.8\linewidth]{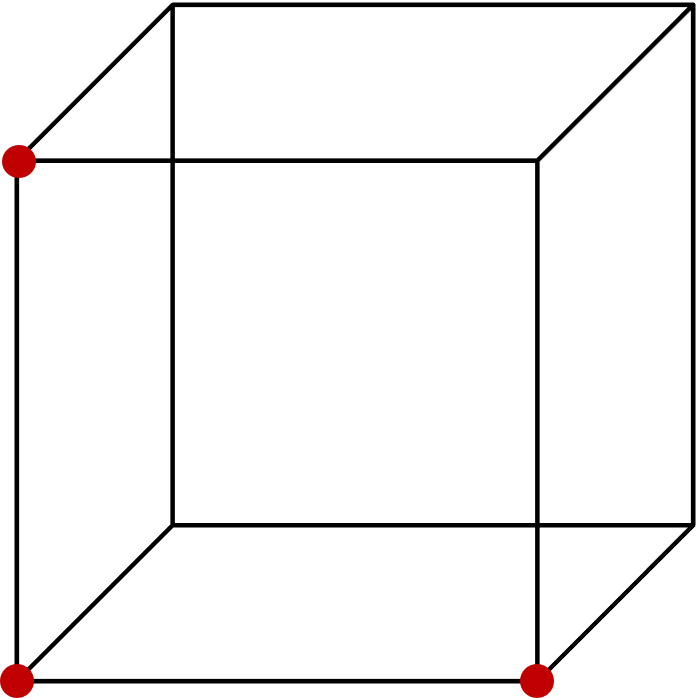}
		\caption{The primal feasible region \\ and three vertices}
		\label{cuttingplane1}
	\end{subfigure}\hspace{0.5cm}
	\centering
	\begin{subfigure}{0.35\linewidth}
		\centering
		\includegraphics[width=0.8\linewidth]{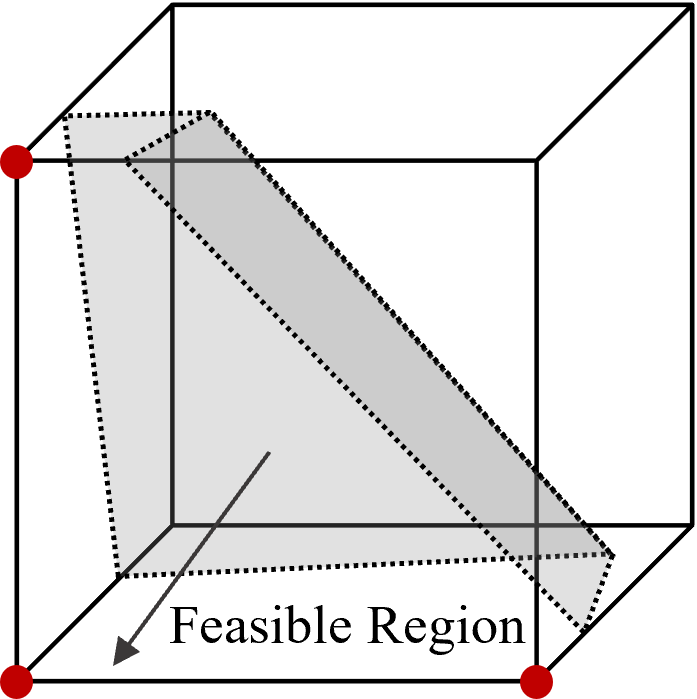}
		\caption{The cutting planes  generated by AE4BV}
		\label{cuttingplane2}
	\end{subfigure}
	\caption{An illustrative example of cutting planes generated by AE4BV. }
	\label{cuttingplane}
\end{figure}

\textbf{Illustrative Example:} To illustrate the effectiveness and mechanism of the cutting plane constraints, we present an example using synthetic data.
We consider a simple dataset $\mathcal{D}_{\boldsymbol{u}} = \{\boldsymbol{u}_1^* = [0 ~ 0 ~ 0]^\top, \boldsymbol{u}_2^* = [0 ~ 1 ~ 0]^\top, \boldsymbol{u}_3^* = [0 ~ 0 ~ 1]^\top\}$ with $p = 3$. As $\boldsymbol{u} \in \mathbb{B}^3$, all data points can be represented as the vertices of a unit cube, as depicted in \picref{cuttingplane1}.  
All data points of the vertices on the unit cube can be regarded as the primal feasible region of a three-dimensional binary vector. 
Assumed that an AE4BV can be trained on the dataset $\mathcal{D}_{\boldsymbol{u}}$ and a polytope $\mathcal{P}_{\boldsymbol{u}}$ can be attained from the decoder parameters $\{\boldsymbol{W}, \boldsymbol{a}\}$, as shown in
\picref{cuttingplane2}. These cutting planes effectively constrain the region to include only the three input vertices while excluding the remaining data points not belonging to $\mathcal{D}_{\boldsymbol{u}}$. By incorporating these cutting planes into MIPs, the feasible region can be effectively shrunk, which facilitates the solution of MIPs by excluding solutions that are less likely to be optimal.

\subsection{Accelerating solving MIPs}\label{sec: accelerating solving}

\begin{figure*}[htbp]
    \centering
    \includegraphics[width=\linewidth]{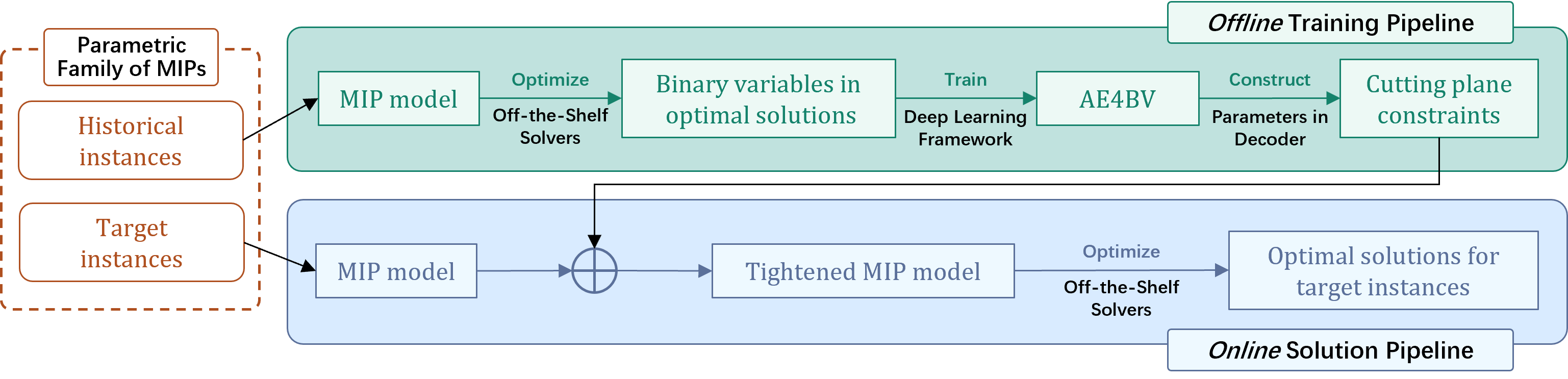}
    \caption{The framework of the approach to accelerating solving MIP problems with AE4BV.}\label{fig: method}
\end{figure*}

Now we arrive at a new L2O framework for accelerating the solution of a parametric family of MIPs, as depicted in \picref{fig: method}. It consists of the offline training pipeline and the online solution pipeline.

\begin{itemize}
    \item \textbf{Offline training}:
    The offline training pipeline involves training AE4BV with data and constructing tightening constraints.   
    First, the optimal solutions of historical instances from a parametric family of MIPs are collected into a training dataset, $\mathcal{D}_{\boldsymbol{u}}$. 
    The AE4BV is then trained on $\mathcal{D}_{\boldsymbol{u}}$, with the network structure proposed in \secref{sec: AE4BV} and the loss function \eqref{eq: MLC}. 
    Once the AE4BV model is trained, the polytope $\mathcal{P}_{\boldsymbol{u}}$ can be constructed based on \eqref{decoder piecewise inequality with Big-M}. 
    The left-hand coefficients in \eqref{decoder piecewise inequality with Big-M} are directly derived from the decoder parameters $\{\boldsymbol{W}, \boldsymbol{a}\}$ and the constant $M$ is given by \eqref{eq: estimation M}. 
    \item \textbf{Online solution}: 
    The online solution pipeline utilizes the cutting plane constraints in $\mathcal{P}_{\boldsymbol{u}}$  to accelerate solving new target MIP problems within the same parametric family. 
    When addressing a new target instance MIP, the cutting planes in  are incorporated into the primal MIP model MIP[$\boldsymbol{\theta}_{\rm new}$], attaining a tightened MIP model T-MIP[$\boldsymbol{\theta}_{\rm new}$]:
\begin{equation}
    \begin{aligned}
        \min_{\boldsymbol{x}, \boldsymbol{u}, \boldsymbol{h}}\quad  & g(\boldsymbol{x}, \boldsymbol{u}; \boldsymbol{\theta}_{\rm new})\\
        \mbox{s.t.}\quad &f(\boldsymbol{x}, \boldsymbol{u}; \boldsymbol{\theta}_{\rm new}) \geq \boldsymbol{0}\\
        &\boldsymbol{W}^\top\boldsymbol{h} + \boldsymbol{a} \geq M\cdot(\boldsymbol{u} - \mathbf{1})\\
        &\boldsymbol{W}^\top\boldsymbol{h} + \boldsymbol{a} \leq M\cdot\boldsymbol{u} \\
        & \boldsymbol{x} \in \mathbb{R} ^ {m}, \boldsymbol{u}\in \mathbb{B} ^{p}, \boldsymbol{h} \in  \mathbb{R}^{d}\\
    \end{aligned}
    \label{mip with constraints}
\end{equation}
    which is still an MIP and can be solved using off-the-shelf solvers. As compared to MIP[$\boldsymbol{\theta}$], the T-MIP[$\boldsymbol{\theta}$] includes additional continuous variables $\boldsymbol{h}$, which correspond to the latent feature in AE4BV. Although this seems to increase the complexity of MIP, the added cuts effectively reduce the feasible region of binary variables and thus help to alleviate the computational burden. 
\end{itemize}

In our approach, the training of AE prepays in part the cost of solving new target problems by encoding the latent structure of binary variables into data-based cutting planes.
As to be empirically shown in Section \ref{sec: 4}, our approach enhances the utilization of heuristics, reduces the number of explored nodes, and, most importantly, significantly reduces the solution time of off-the-shelf MILP solvers. The price we have to pay is only a slight loss of optimality due to a tightened problem to be resolved. 
 
It is worth noting that our approach is distinct from known end-to-end L2O methods that directly predict the binary variables through supervised learning.  
When there are numerous decision variables and complex coupling constraints in MIP, end-to-end learning can struggle to provide feasible solutions, while our method offers an alternative useful option for promoting the overall efficiency of solving MIPs. There are also some relevant works that handle complex constraints in MIPs in a data-driven fashion. For instance,
\cite{xavier2021learning} devised a problem-specific heuristic to construct hand-crafted cutting planes to tighten SCUC problems. 
In \cite{pineda2020data}, another data-driven strategy was proposed to remove some constraints in the primal problem to obtain a more relaxed model. Our method is in a similar spirit to \cite{xavier2021learning} but offers a more general pipeline that can automatically learn tightening constraints from data to accelerate the solution process of MIP solvers. In this sense, our method enjoys full generality in systematically integrating unsupervised learning and operations research, thereby expanding the scope of the current L2O paradigm. 

\section{Computational Experiment}
\label{sec: 4}
\hbadness=10000

In this section, we conduct a comprehensive analysis of the performance of our proposed method applied to a batch process scheduling problem.
The experiments are carried out on a desktop computer with an Intel(R) Core(TM) i7-12700F @ 2.10 GHz and an NVIDIA GeForce RTX 2060 graphics card, operating under a Windows 10 environment. 
All algorithms are implemented in Python 3.10.4 and the codes are made available at \href{https://github.com/qushiyuan/AE4BV}{https://github.com/qushiyuan/AE4BV}.

\subsection{Problem description}
\label{results:bps}
\begin{figure}[htbp]
    \centering
    \includegraphics[width=\linewidth]{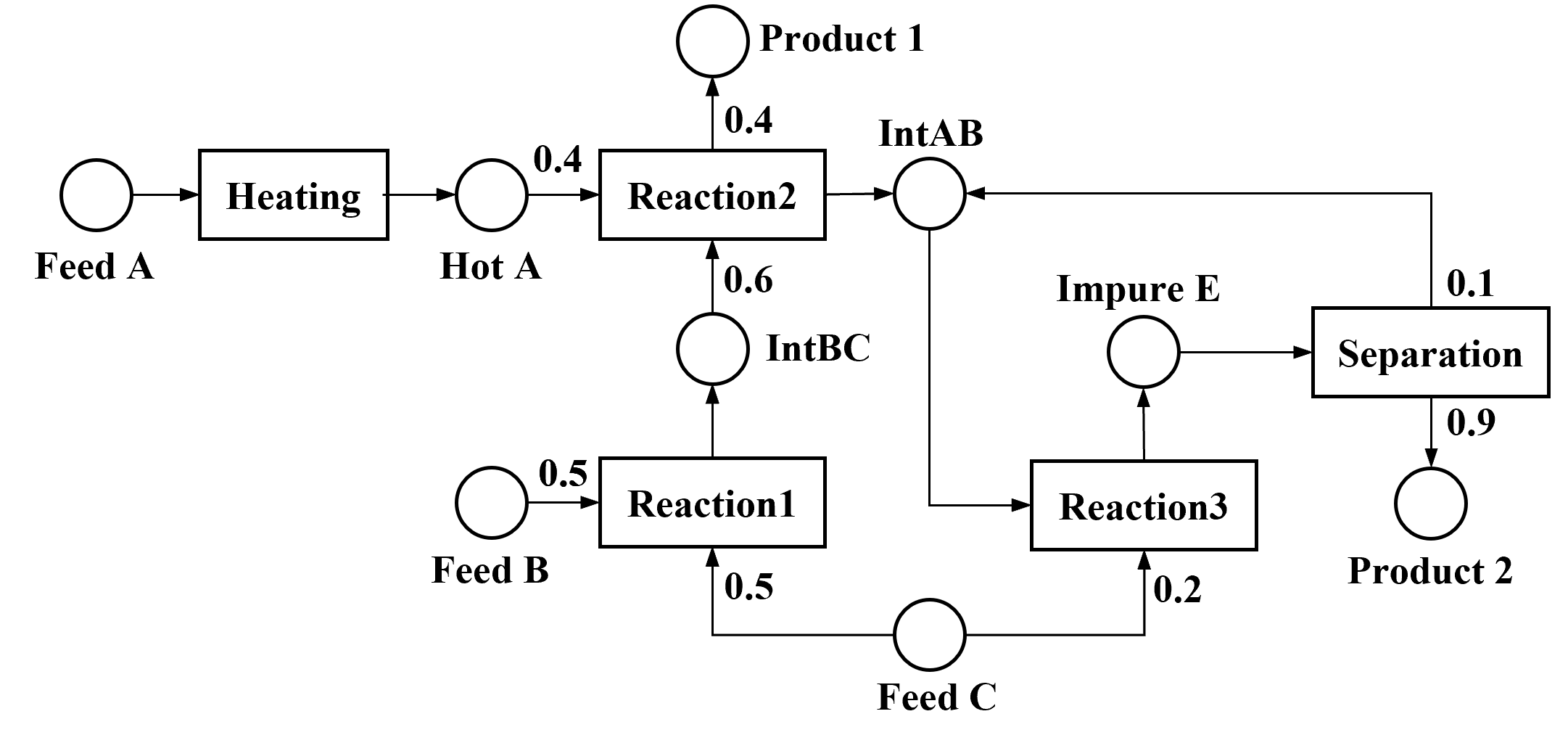}
    \caption{The STN representation of the batch process problem. }\label{fig:bps_scheduling}
\end{figure}

Over the past few decades, the batch process scheduling problem has been recognized as one of the most important challenges in industrial scheduling operations.
Formulated as MILPs, these scheduling problems need to be resolved repeatedly. In our experiments, we adopt the continuous-time formulation of batch process scheduling problems proposed by \cite{ierapetritou1998effective}.  
The problem focuses on the production of two products through five processing tasks, including heating, three reactions (denoted as reactions 1, 2, and 3), and separation. The state-task network (STN) representation of this batch process is shown in \picref{fig:bps_scheduling}.
It involves 9 distinct material stages and 8 units, each being capable of performing specific tasks within the process. 

The entire batch process scheduling problem can be cast as a complicated MILP having the form of \eqref{standard mip}. The full problem formulation, along with the relevant data, can be found in \ref{appendix:bps}, where the objective is to maximize the total sales of products, and the binary variables indicate the starting of each task and the units utilized in the process. \tabref{tab:scale of MIP model} presents detailed information on the size of the induced MILP model, including the number of rows, columns, nonzero elements, and the number of continuous and binary variables. 

\begin{table*}[htbp]
  \caption{The size of the primal MILP model and tightened MILP model ($N = 1000, \mathcal{E} = 5\%$) of batch processing. The number of cutting planes, explored nodes, and solution time is recorded based on the solution process using Gurobi. 
  The better performances are highlighted in \underline{underline}.}
  \centering
  \adjustbox{width=\linewidth}{
  
    \begin{tabular}{lccccccccccc}
    \toprule
          & \multicolumn{4}{c}{Primal model}      & \multicolumn{4}{c}{Presolved model}   & \multirow{2}[4]{*}{\makecell[c]{~Cutting~ \\planes}} & \multirow{2}[4]{*}{\makecell[c]{~Explored~ \\nodes}} & \multirow{2}[4]{*}{\makecell[c]{~Solution~ \\time}}\\
\cmidrule(lr){2-5}  \cmidrule(lr){6-9}        & ~Rows~  & ~Columns~ & ~Nonzeros~ & ~Binaries~ & ~Rows~  & ~Columns~ & ~Nonzeros~ & ~Binaries~ &       &   & \\
    \midrule
    Primal MILP & 3431 & 3834 & 13588 & 216  & 2317 & 565   & 9932 & 122   & 405   & 127282 & 33.47s\\
    ~Tightened MILP~ & 3863  & 3844  & 18340  & 216   & 2423  & 554 & 12137 & 110 & \underline{212} & \underline{1621} & \underline{2.25s}\\
    \bottomrule
    \end{tabular}}%
  \label{tab:scale of MIP model}%
\end{table*}%

\subsection{Offline AE4BV training}
For L2O, we generate a training dataset $\mathcal{D}_{\boldsymbol{u}}$ with $N=1000$ samples by imposing perturbations on problem parameters of batch process scheduling and then solving the induced problem instances to optimality using Gurobi 10.0.0. We consider three types of parameter perturbations: variation in the constant term $T_{ij}^c$ of task processing time, variation in the variable term $T_{ij}^v$ of processing time per unit of material for tasks, and variation in the unit price $PRICE_s$ of products. Three different levels of perturbation, i.e. $\mathcal{E} =  5\%, 10\%, 20\%$, are considered, which quantify how much the relevant parameters are varied. Concretely, each parameter is varied by multiplying its nominal value with a random factor within [0.95, 1.05], [0.90, 1.10], [0.80, 1.20], respectively, corresponding to $\mathcal{E} =  5\%, 10\%, 20\%$. 

Before proceeding with AE4BV, we first investigate the performance of existing end-to-end L2O methods. Specifically, we choose the optimal classification tree (OCT) \citep{bertsimas2017optimal}, which is trained on a dataset $\mathcal{D}_{\boldsymbol{u}}$ with $N=1000$ and $\mathcal{E} = 5\%$, and then used to straightforwardly predict the optimal solution to the scheduling problem. 
Its performance is evaluated on another test set including $N_{\rm test} = 200$ instances; 
however, it turns out that all predicted solutions are infeasible. This highlights the critical challenge of directly predicting optimal solutions for process scheduling problems, as well as the need for MILP solvers to ensure the feasibility and optimality of solutions.  

Next, we train AE4BV to accelerate the solution of MILP solvers. The hyperparameters of AE4BV are chosen as follows: the number of hidden layers ([6, 12]), the number of elements of the feature vector $\boldsymbol{h}$ ([10, 30]), the learning rate ([$10^{-4}$, $10^{-1}$]), and the number of training epochs ([100, 2000]). 
The final configuration of these hyperparameters is presented in \tabref{table:hyperparameter of AE4BV in BPS}. To evaluate the performance of AE4BV, we conduct experiments on a test dataset comprising $N_{\rm test} = 200$ instances for each choice of $\mathcal{E}$. 
The following metrics are used: 
\begin{itemize}
\item\textit{Hamming loss} (HL) \citep{tsoumakas2008multi} is a widely used performance metric in MLC:   
\begin{equation}
    \text{HL} = \frac{1}{N_{\rm test}}\sum_{n = 1}^{N_{\rm test}} \frac{\Vert \boldsymbol{u}^*_n \oplus \hat{\boldsymbol{u}}^*_n \Vert_1}{P}
    \label{HL}
\end{equation}
where $ \oplus $ denotes the element-wise XOR operation and $\hat{\boldsymbol{u}}^*_n$ is the reconstructed input of AE4BV with $ \boldsymbol{u}^*_n$ used as the input. 
The smaller the Hamming loss, the better the performance of feature extraction and input reconstruction.

\item\textit{Probability of preserving optimality} (PPO) quantifies the probability of whether the polytope \eqref{decoder piecewise inequality with Big-M} built from the learned AE4BV manages to enclose the optimal binary variables $\{ \boldsymbol{u}^*_n \}$: 
\begin{equation}
    \text{PPO} = \frac{1}{N_{\rm test}}\sum_{n = 1}^{N_{\rm test}} \mathbb{I}(\boldsymbol{u}^*_n \in \mathcal{P}_{\boldsymbol{u}})
    \label{HL}
\end{equation}
where $\mathbb{I}(\cdot)$ is the indicator function.
If the optimal solution is contained within the polytope $\mathcal{P}_{\boldsymbol{u}}$, optimality can be ensured when solving the tightened MIP. Otherwise, there can be a certain loss of optimality.

\end{itemize}

\begin{table}[h!]
    \caption{Hyperparameter configuration of the AE4BV in the experiment on batch processing. }
    \label{table:hyperparameter of AE4BV in BPS}
    \centering
    \adjustbox{width=\linewidth}{
    \begin{tabular}{lc}
    \toprule
    Hyperparameter & Values / Type\\
    \midrule
        Epochs & 500\\
        Optimizer & Adam\\
        Activate function & LeakyReLU\\
        Learning rate & $2\times 10^{-4}$ \\
        Dropout & 0.2 \\
        Number of elements in $\boldsymbol{h}$ & 20 \\
        Neurons of encoder layers~~~~ & ~~~~[20, 216, 40, 216, 120, 216, 180, 216]~~~~\\
    \bottomrule
    \end{tabular}}
\end{table}
\begin{table}[h!]
  \centering
  \caption{Performance of AE4BV on batch processing.}
  \adjustbox{width=\linewidth}{
    \begin{tabular}{lcccccc}
    \toprule
          & \multicolumn{2}{c}{$\mathcal{E} = 5\%$} & \multicolumn{2}{c}{$\mathcal{E} = 10\%$} & \multicolumn{2}{c}{$\mathcal{E} = 20\%$} \\
\cmidrule(lr){2-3} \cmidrule(lr){4-5} \cmidrule(lr){6-7}  $N$ & ~HL (\%) ~& ~PPO (\%)~ & ~HL (\%)~ & ~PPO (\%)~  & ~HL (\%)~ & ~PPO (\%)~\\
    \midrule
    $200$ & 3.16  & 95    & 3.68  & 85.5  & 5.47  & 75.5 \\
    $500$ & 1.22  & 97.5  & 1.66  & 96    & 3.14  & 86 \\
    $1000$ & 0.39  & 98    & 0.8   & 97    & 1.94  & 92 \\
    \bottomrule
    \end{tabular}}%
  \label{tab:performance of AE4BV}%
\end{table}%
\begin{table}[h!]
  \centering
  \caption{Performance of AE4BV without skip connections on batch processing. }
  \adjustbox{width=\linewidth}{
    \begin{tabular}{lcccccc}
    \toprule
          & \multicolumn{2}{c}{$\mathcal{E} = 5\%$} & \multicolumn{2}{c}{$\mathcal{E} = 10\%$} & \multicolumn{2}{c}{$\mathcal{E} = 20\%$} \\
\cmidrule(lr){2-3} \cmidrule(lr){4-5} \cmidrule(lr){6-7} $N$ & ~HL (\%)~ & ~PPO (\%)~ & ~HL (\%)~ & ~PPO (\%)~  & ~HL (\%)~ & ~PPO (\%)~\\
\midrule
    $200$ & 5.23  & 89.5  & 6.19  & 73    & 8.31  & 44.5 \\
    $500$ & 3.43  & 90.5  & 4.55  & 87    & 6.05  & 49.5 \\
    $1000$ & 2.48  & 91 & 3.28  & 75    & 4.48  & 43 \\
    \bottomrule
    \end{tabular}}%
  \label{tab:performance of AE4BV without skip connections}%
\end{table}%

The overall performance of AE4BV on the test dataset, under different training set sizes $N$ and $\mathcal{E}$, is summarized in \tabref{tab:performance of AE4BV}. 
It shows the desirable generalization capability of AE4BV for binary data. 
Specifically, the results indicate that AE4BV performs better with larger $N$ and smaller $\mathcal{E}$. 
Consequently, it is advisable to collect as many data samples as possible to enhance the effectiveness of AE4BV.  
As $\mathcal{E}$ becomes larger, both HL and PPO indices become worse. 
This implies using a larger $\mathcal{E}$ causes a more scattering distribution of binary variables, resulting in worse performance of AE4BV in dimension reduction. 
Besides, when $N$ is sufficiently large, PPO remains high albeit still less than 100\%. 
As to be shown later, the optimality gap between the tightened MILP and the primal problem is rather small. 

Next, we perform an ablation study to shed light on the impact of skip connections on the performance of AE4BV. 
The results of AE4BV without skip connections are presented in \tabref{tab:performance of AE4BV without skip connections}. 
Compared to the full AE4BV, the performance of AE4BV without skip connections degrades significantly. 
In particular, when either $N$ or $\mathcal{E}$ becomes larger, the performance gap between the two networks becomes wider.
These results highlight the pivotal role of skip connections in ensuring the desirable reconstruction performance of AE4BV for binary input data.

\subsection{Online solution performance}\label{sec:online solving performance bps}
Next, we study the optimization performance of the proposed unsupervised learning scheme on the batch process scheduling problem. 
After training AE4BV, we obtain data-based cutting plane constraints and use them to construct tightened MILP models for new problem instances. 

In this subsection, we train AE4BV with $N = 1000$ and different values of $\mathcal{E}$. Then the tightened MILPs are constructed by adding cutting planes to a set of $N_{\rm test} = 200$ test instances of the batch process scheduling problem, which are then solved using various off-the-shelf MILP solvers, including Gurobi, COPT, and SCIP, each configured with its own default setting.  
Considering the performance discrepancies among the three solvers, the solution time for Gurobi and COPT is computed as the time required to attain global optimality, while for SCIP the solution time is defined as the duration needed for reaching an optimality gap of 0.02.

\begin{table*}[htbp]
  \centering
  \caption{Average, maximum, and standard deviation of solution time using different solvers with primal and tightened MILPs. }
  \adjustbox{width=\linewidth}{
    \begin{tabular}{llccccccccc}
    \toprule
    \multicolumn{2}{c}{\multirow{2}[3]{*}{}} & \multicolumn{3}{c}{Gurobi} & \multicolumn{3}{c}{COPT} & \multicolumn{3}{c}{SCIP} \\
    \cmidrule(lr){3-5}\cmidrule(lr){6-8}\cmidrule(lr){9-11}    \multicolumn{2}{c}{} & $t_{\text{avg}}/s$ & $t_{\text{max}}/s$ & $t_{\text{std}}/s$  & $t_{\text{avg}}/s$ & $t_{\text{max}}/s$ & $t_{\text{std}}/s$ & $t_{\text{avg}}/s$ & $t_{\text{max}}/s$ & $t_{\text{std}}/s$  \\
    \midrule
    \multirow{3}[0]{*}{$\mathcal{E} = 5\%$} ~~~~~~~~& Primal MILP~~~~~~~~~ & 20.47 & 154.16 & 24.84 & 50.80  & 467.49 & 59.49 & 73.73 & 709   & 96.30 \\
          & Tightened MILP & 2.10 & 4.41 & 0.47 & 3.13 & 5.19 & 0.67 & 12.29 & 31 & 3.87\\
          & Speedup (\%) & ~~~89.76~~~ & ~~~97.14~~~ & ~~~98.10~~~ & ~~~93.83~~~ &  ~~~98.89~~~ & ~~~98.87~~~ & ~~~83.33~~~ & ~~~95.63~~~ & ~~~95.98~~~ \\
    \midrule
    \multirow{3}[0]{*}{$\mathcal{E} = 10\%$} & Primal MILP & 27.16 & 205.74 & 36.55 & 86.29 & 932.00   & 130.29 & 184.42 & 2151  & 362.15 \\
          & Tightened MILP & 6.76 & 17.82 & 3.03 & 9.76 & 32.35 & 6.39 & 31.92 & 108 & 21.11 \\
          & Speedup (\%) & 75.13 & 91.34 & 91.70 & 88.69 &  96.53 & 95.10 & 82.69 & 94.98 & 94.17 \\
    \midrule
    \multirow{3}[0]{*}{$\mathcal{E} = 20\%$} & Primal MILP & 34.59 & 384.10 & 54.19 & 111.72 & 1027.87 & 159.82 & 407.51 & 5057  & 853.11 \\
          & Tightened MILP & 12.87 & 60.35 & 6.94 & 30.35 & 83.66 & 19.45 & 83.04 & 268 & 61.05 \\
          & Speedup (\%) & 62.79 & 84.29 & 87.20 & 72.83 & 91.86 & 87.83 & 79.62 & 94.70 & 92.84 \\
    \bottomrule
    \end{tabular}}%
  \label{tab:solving time}
\end{table*}%

\begin{figure}[h!]
    \begin{subfigure}{\linewidth}
        \centering
        \includegraphics[width=\linewidth]{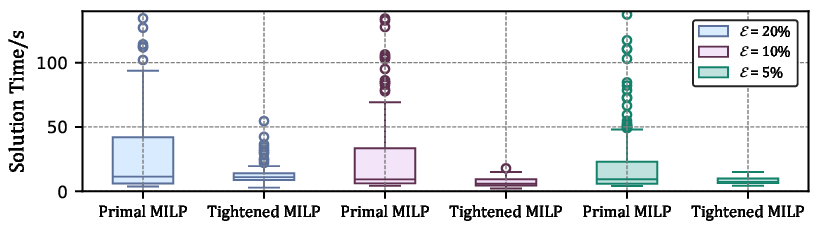}
        \caption{The solution time using Gurobi}
        \label{fig:exp_time_gurobi}
    \end{subfigure}
    \begin{subfigure}{\linewidth}
        \centering
        \includegraphics[width=\linewidth]{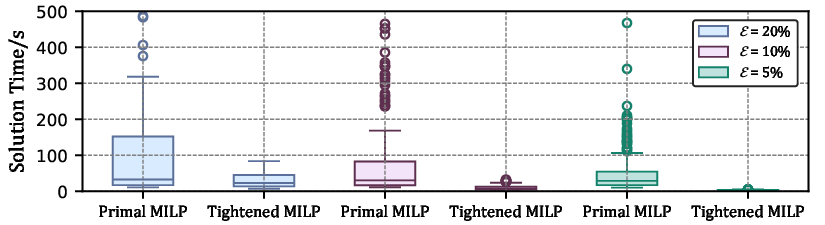}
        \caption{The solution time using COPT}
        \label{fig:exp_time_copt}
    \end{subfigure}
    \begin{subfigure}{\linewidth}
        \centering
        \includegraphics[width=\linewidth]{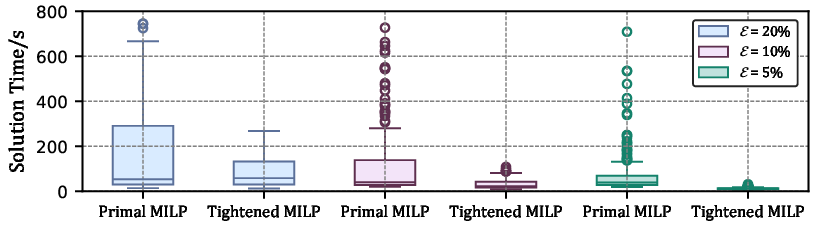}
        \caption{The solution time using SCIP with MIP gap of 0.02}
        \label{fig:exp_time_scip}
    \end{subfigure}
    \caption{The solution time using three solvers under different problem settings of batch processing. }
    \label{fig:exp_time}
\end{figure}

\picref{fig:exp_time} presents the boxplots of solution time under different problem settings with primal and tightened MILPs through different solvers. 
At first glance, the time required to solve tightened MILPs is significantly less than that of primal MILPs in all cases. 
This indicates that our proposal is a general and effective approach to accelerating the solution of a parametric family of MIPs, thanks to the reduction of the feasible region through unsupervised learning. 
Additionally, the solution time of tightened MILPs exhibits much lower variations and lighter tails than that of primal MILPs. This indicates that the proposed unsupervised learning approach is particularly helpful for managing the solution cost of a family of MILPs and avoiding some extremely slow convergence of off-the-shelf solvers. 

\tabref{tab:solving time} presents the detailed statistics containing the average solution time $t_{\rm avg}$, maximum solution time $t_{\rm max}$, and the standard deviation of the solution time $t_{\rm std}$.
Additionally, the percentage of speedup achieved by the tightened MILP models is also shown. 
Across all cases, the improvement achieved by our approach is approximately an order of magnitude. 

\begin{figure}[htbp]
    \centering
    \includegraphics[width = \linewidth]{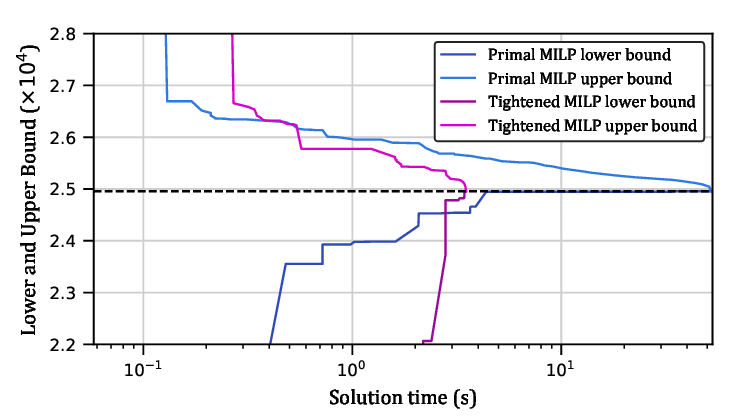}
    \caption{The solution process of batch processing using COPT. The black dashed line indicates the optimal objective value.}\label{fig:copt solution process}
\end{figure}

To gain more insight into the effectiveness of our method, we further investigate a particular problem instance whose solution time is close to the average.
The relevant information of the tightened MILP model and its solution process using Gurobi, is presented in \tabref{tab:scale of MIP model}. 
It can be seen that the size of the tightened MILP does not increase significantly as compared to the primal MILP. Interestingly, after pre-solving with Gurobi, the number of binary variables of the tightened MILP is even smaller than that of the primal MILP. 
The number of generated cutting planes and the nodes explored throughout the B\&C process are also presented. 
Compared to the primal MILP model, solving the tightened MILP model involves much fewer generated cutting planes and explored nodes in B\&C, and most notably, a significantly reduced solution time. 
The fewer generated cutting planes and the fewer nodes explored demonstrate that the complexity of solving the tightened MILP gets considerably lowered, which eventually leads to the reduction of solution time. 

We further use COPT solver to solve both the tightened problem and the primal problem for a particular instance, and plot the profiles of their upper and lower bounds during solution processes in \picref{fig:copt solution process}, where a logarithmic time scale is used for better clarity. Recall that the problem we solve is a maximization problem. Obviously, when solving the primal MILP, the primary limiting factor is the slow convergence of the upper bound, which stems from continually solving relaxed problems during B\&C process.  
By adding polytopic constraints, the tightened MILP shows a much faster decrease in the upper bound. This implies that the integration of tightening constraints effectively helps to improve the quality of relaxations during B\&C process. Besides, the lower bound also increases faster when solving the tightened MIP, which is mainly due to a smaller space for finding feasible solutions. 
This highlights the effectiveness of our approach in enhancing the efficiency of B\&C process in MILP solvers.

\begin{table}[h!]
  \centering
  \caption{The gap between the optimal objective values of the primal MILP and the tightened MILP under different problem settings of batch processing. } 
  \adjustbox{width=\linewidth}{
    \begin{tabular}{lcccccc}
    \toprule
          \multicolumn{1}{l}{\multirow{2}[3]{*}{$\mathcal{E}$}~~~~~~~~} & {\multirow{2}[3]{*}{~~~~~$\overline{\text{Gap}} (\%)$~~~~~}} & \multicolumn{5}{c}{~~~~~~~Percentage within Gap (\%)~~~~~~} \\
\cmidrule(lr){3-7}    &       & \multicolumn{1}{c}{~~~$\leq$1\%~~} & \multicolumn{1}{c}{~~$\leq$2\%~~} & \multicolumn{1}{c}{~~$\leq$3\%~~} & \multicolumn{1}{c}{~~$\leq$4\%~~} & \multicolumn{1}{c}{~~$\leq$5\%~~~} \\
    \midrule
    $5\%$   & 0.09  & ~100~   & ~100~   & ~100~   & ~100~   & ~100~ \\
    $10\%$  & 0.42  & 85.5  & 99.5  & 100   & 100   & 100 \\
    $20\%$  & 1.04  & 60    & 76.5  & 94    & 99    & 100 \\
    \bottomrule
    \end{tabular}}
  \label{tab:optimality keeping}%
\end{table}%

Finally, we investigate the solution quality of tighted problems. 
For each $\mathcal{E}$, \tabref{tab:optimality keeping} shows the average gap between optimal objective values of the tightened MILP model and the primal MILP model, as well as the percentage of different levels of optimality achieved. 
When $\mathcal{E} = 5\%$, the loss of  optimality resulting from the tightened MILP model is at most 0.68\%, which is rather minor.
As $\mathcal{E}$ increases, the quality of solutions moderately decreases.
When $\mathcal{E} = 20\%$, our method only induces an average optimality loss of about 1\%. To summarize, even if the proposed approach does not guarantee the optimality, it significantly reduces the expensive computational cost while maintaining a fairly high quality of solutions to the rightened problem. 

\section{Conclusion}

In this work, we presented a novel L2O scheme that lies at the intersection of unsupervised learning and operations research. It makes use of AE4BV to disentangle the latent structure of binary variables of optimal solutions and accelerate the solution process of a parametric family of MIPs. By adding a class of cutting plane constraints generated from a trained AE4BV, tightened MIPs can be formulated for new problem instances and then resolved by off-the-shelf solvers with a higher efficiency. 
Different substantially from the existing machine learning methods for optimization, our approach offers a new perspective in L2O and paves the way for integrating data-driven insights into optimization frameworks.  
When applied to a benchmark batch process scheduling problem, our approach achieved a reduction in solution time by an order of magnitude while achieving an optimality loss of less than 1\%.
These findings highlight the effectiveness of AE4BV in enhancing the efficiency and practicality of solving MIPs. 

For future research, we will explore how to improve the generalization capability of AE4BV and reduce the loss of optimality of tightened problems. The usage of our approach in tackling more general MIP problems is also worth investigating, including mixed integer convex programming and mixed integer nonlinear programming. 

\section*{Acknowledgments}

This work was supported by the National Natural Science
Foundation of China under Grants 62373211 and 62327807.

\appendix

\section{Batch process problem formulation}
\hbadness=10000
\label{appendix:bps}
Here we present the full MILP formulation of the batch process problem which we use in the computational experiments in \secref{results:bps}. 
This batch process is composed of $8$ units, $5$ tasks corresponding to fixed units, and $9$ different material stages required for the production of $2$ final products. 
The data for this problem are presented in \tabref{table:data for bps unit-task} and \tabref{table:data for bps state}.

\begin{table}[h!]
    \caption{Unit-task data for the batch process problem}
    \label{table:data for bps unit-task}
    \centering
    \adjustbox{width=\linewidth}{
    \begin{tabular}{lccc}
    \toprule
    Unit & Capacity& Suitability&Mean processing time\\
    \midrule
        Heater 1 & 100 & Heating & 4.0 \\
        Heater 2 & 120 & Heating & 4.0 \\
        Heactor 1~~ & ~70~ & ~~~~Reaction 1, 2, 3~~~~& 4.0 \\
        Reactor 2 & 80 & Reaction 1, 2, 3 & 3.0 \\
        Reactor 3 & 70 & Reaction 1, 2, 3 & 4.0 \\
        Reactor 4 & 120 & Reaction 1, 2, 3 & 5.0 \\
        Still 1 & 200 & Separation & 5.0 \\
        Still 2 & 150 & Separation & 5.0 \\
    \bottomrule
    \end{tabular}
    }
\end{table}
\begin{table}[h!]
    \caption{Material state data for the batch process problem}
    \label{table:data for bps state}
    \centering
    \adjustbox{width=\linewidth}{
    \begin{tabular}{lccc}
    \toprule
    State &Storage capacity&~~~~~~Initial amount~~~~~~&~~~~Price~~~~~~~~\\
    \midrule
        Feed A & Unlimited & 1000 & 0 \\
        Feed B & Unlimited & 800 & 0 \\
        Feed C & Unlimited & 800 & 0 \\
        Hot A & 100 & 0 & 0 \\
        IntAB & 200 & 0 & 0 \\
        IntBC & 150 & 0 & 0 \\
        Impure E & 200 & 0 & 0 \\
        Product 1~~~~~ & ~~~~~Unlimited~~~~~~ & 0 & 25 \\
        Product 2 & Unlimited &0 & 30 \\
    \bottomrule
    \end{tabular}}
\end{table}

The MILP formulation is based on the continuous-time representation with $n \in N$ event time points. Furthermore, let $i \in I, j \in J, s \in S$ represent the tasks, units, and material states respectively. 
Let $I_j$ be the tasks that can be performed in unit $j$, and $J_i$ be the units that are suitable for performing task $i$. 
Let $Is$ be the tasks that either produce or consume state $s$. 
The following constant parameters will be used in MILP formulation:

\begin{table}[ht]
\centering
\adjustbox{width=\linewidth}{
\begin{tabular}{cl}
\toprule
\textbf{Parameters} & \textbf{Description} \\ \midrule
$C_{ij}^{\rm max}$ & Maximum capacity of unit $j$ processing task $i$ \\ 
$S_s^{\rm max}$ & Maximum storage capacity for state $s$ \\
$P_{si}^p$ & Proportion of state $s$ produced in task $i$ \\ 
$P_{si}^c$ & Proportion of state $s$ consumed in task $i$\\
$T_{ij}^c$ & Constant term of processing time of task $i$ at unit $j$ \\ 
$T_{ij}^v$ & Variable term of processing time per material unit of task $i$ at unit $j$ \\ 
$A_s$ & The initial amount of state $s$ \\ 
$H$ & The time horizon \\ 
$PRICE_s$ & Unit price of state $s$ \\ \bottomrule
\end{tabular}}
\caption{List of Variables and Descriptions for MILP model of the batch process problem}
\end{table}

For each event point $n \in N$, various variables are defined to describe the task, unit, and state. 
The main decision variables for this problem are binary variables, $m(i, n)$ and $u(j, n)$. 
$m(i, n)\in \{0, 1\}$ indicates whether the task $i$ is started at the event point $n$, and $u(j, n)\in \{0, 1\}$ indicates whether unit $j$ is utilized at event point $n$. 
Other auxiliary continuous variables are used to describe the amount of state and real-time event points. 
$p(i, j, n)$ indicates the total amount of material undertaking task $i$ in unit $j$ at event point $n$.  
$s(s, n)$ indicates the sold amount if state $s$ at event point $n$.
$sa(s, n)$ indicates the total amount of state $s$ at event point $n$. 
Finally, $t^s(i, j, n)$ and $t^f(i, j, n)$ are respectively the start time and finish time of task $i$ in unit $j$ at event point $n$. 
Given the parameters and variables above, the batch process problem can be formulated as a MILP problem:
\begin{flalign}
    {\rm max\ } & \sum_{s\in S}\sum_{n\in N}{PRICE_s \cdot s(s,n)} & \label{bpsc:obj}\\
    {\rm s.t.} & \sum_{i\in{I_j}}{m(i,n)}=u(j,n), & \forall {j\in{J}},n\in{N}, \label{bpsc:allocation}
\end{flalign}
\begin{flalign}
    \qquad & sa(s, 0)=A_s, & \forall{s\in{S}}, \label{bpsc:init}\\
    \qquad & p(i,j,n)\leq{sa(s,n)}, &\forall{i\in{I}},j\in{J_i},\notag \\
        & & \forall n\in{N},s\in{I_s}, \label{bpsc:quantity}\\
    \qquad & p(i,j,n)\leq{C_{ij}^{\rm max}}m(i,n), &     \forall{i\in{I}},\notag\\
        & & \forall j\in{J_i},n\in{N}, \label{bpsc:capacity}\\
    \qquad & sa(s,n)\leq{S_s^{\rm max}}, &\forall{s\in{S}},n\in{N}, \label{bpsc:storage}\\
    \qquad & sa(s,n)=sa(s,n-1) \notag\\
        & \hspace{1em} -s(s,n) \notag \\
        & \hspace{1em} + \sum_{i\in{I_s}}{P_{si}^p}\sum_{j\in{J_i}}{p(i,j,n-1)} \notag\\
        & \hspace{1em} + \sum_{i\in{I_s}}{P_{si}^c}\sum_{j\in{J_i}}{p(i,j,n)},&  \forall{s\in{S}},n\in{N}, \label{bpsc:material}\\
    \qquad & t^f(i,j,n)=t^s(i,j,n)  \notag \\
        & \hspace{1em} +T_{ij}^cm(i,n)+T_{ij}^v,\notag \\
        & & \forall{i\in{I}},\notag \\
        & & \forall j\in{J_i},n\in{N}, \label{bpsc:reaction time}\\
    \qquad & t^s(i,j,n+1) \geq t^f(i,j,n) \notag\\
        & \hspace{1em} - H(2-m(i,n)-u(j,n)),  & \notag\\
        & & \forall{i\in{I}},j\in{J_i}, \notag \\
        & & \forall n\in{N},n\neq{N}, \label{bpsc:reaction seq 1}\\
    \qquad & t^s(i,j,n+1)\geq{t^s(i,j,n)},\notag\\
        & & \forall{i\in{I}},j\in{J_i},\notag \\
        & & \forall n\in{N},n\neq{N}, \label{bpsc:reaction seq 2}\\
    \qquad & t^f(i,j,n+1)\geq{t^f(i,j,n)},\notag\\
        & & \forall{i\in{I}},j\in{J_i},\notag \\
        & & \forall n\in{N},n\neq{N}, \label{bpsc:reaction seq 3}\\
    \qquad & t^s(i,j,n+1)\geq{t^f(i',j,n)} \notag \\
        & \hspace{1em} -H(2-m(i',n)-u(j,n)),\notag \\
        & &\forall{j\in{J}},i\in{I_j},\notag \\
        & & \forall i'\in{I_j},i\neq{i'},\notag \\
        & & \forall n\in{N},n\neq{N}, \label{bpsc:reaction seq 4}\\
    \qquad & t^s(i,j,n+1)\geq{t^f(i',j',n)}\notag \\
        & \hspace{1em} -H(2-m(i',n)-u(j',n)),\notag\\
        & & \forall{j,j'\in{J}},i\in{I_j},\notag \\
        & & \forall i'\in{I_j'},i\neq{i'},\notag \\
        & & \forall n\in{N},n\neq{N}, \label{bpsc:reaction seq 5}\\
    \qquad & t^s(i,j,n+1)\geq\notag \\
        & \hspace{1em} \sum_{n'\in{N},n'\leq{n}}\sum_{i'\in{I_j}}{t^f(i',j,n')},\notag\\
        & \hspace{1em} -\sum_{n'\in{N},n'\leq{n}}\sum_{i'\in{I_j}}{t^s(i',j,n')},\notag\\
        & & \forall{i\in{I}},j\in{J_i},\notag\\
        & & n\in{N},n\neq{N}, \label{bpsc:reaction seq 6}\\
    \qquad & t^f(i,j,n)\leq{H},& \forall{i\in{I}},\notag \\
        & & \forall j\in{J_i},n\in{N}, \label{bpsc:time1}\\
    \qquad & t^s(i,j,n)\leq{H},& \forall{i\in{I}},\notag \\
        & & \forall j\in{J_i},n\in{N}. \label{bpsc:time2}
\end{flalign}

The objective function \eqref{bpsc:obj} considers the sales of the products. 
Equation \eqref{bpsc:allocation} enforces that each unit can only perform one task at every event point. 
Equation \eqref{bpsc:init} enforces the initial amount of each state. 
Equation \eqref{bpsc:quantity}-\eqref{bpsc:storage} enforces the validity and the capacity and storage limits of materials. 
Equation \eqref{bpsc:material} guarantees the conservation of matter throughout processing. 
Equation \eqref{bpsc:reaction time}-\eqref{bpsc:reaction seq 6} enforces the processing time and sequence. 
Equations \eqref{bpsc:time1} and \eqref{bpsc:time2} guarantee the tasks are all performed on the time horizon.

\bibliographystyle{elsarticle-harv} 
\bibliography{elsarticle-num}






\end{document}